\documentclass[12pt]{article}
 \usepackage{latexsym}
 \usepackage{amssymb}
 \usepackage{graphicx}
 \usepackage{amsmath}

 \newtheorem{Theorem}{Theorem}
\newtheorem{Definition}{Definition}
\newtheorem{Proposition}{Proposition}

\newtheorem{Lemma}{Lemma}
\newtheorem{Corollary}{Corollary}

\newtheorem{Example}{Example}

\newcommand{\A}{{\cal A}}
\newcommand{\B}{{\cal B}}
\newcommand{\C}{{\cal C}}

\newcommand{\I}{{\cal I}}
\newcommand{\LL}{{\cal L}}

\newcommand{\x}{{\bf x}}
\newcommand{\y}{{\bf y}}

\newcommand{\0}{{\bf 0}}

\newcommand{\qed}{\nobreak \ifvmode \relax \else
      \ifdim\lastskip<1.5em \hskip-\lastskip
      \hskip1.5em plus0em minus0.5em \fi \nobreak
      \vrule height0.75em width0.5em depth0.25em\fi}

\def \ep{\hbox{ }\hfill$\Box$}

\begin{document}
\title{Further results on $B$-tensors with application to the location of real eigenvalues}

\author{
Lu Ye
\thanks{ College of Economics and Management and Research Center for Ecological Civilization, Zhejiang Sci-Tech University, Hangzhou, P.R. China. Email: zjwzajyl@126.com. This author's work is supported by the Applied Economics Base Fund, Zhejiang Sci-Tech University. } \quad
Zhongming Chen
\thanks{ Corresponding author. Department of Electrical and Electronic Engineering, The University of Hong Kong, Hong Kong.
Email: zmchen@eee.hku.hk.
This author's work was partially done
when he was visiting the Hong Kong Polytechnic University.},
}

\date{\today} \maketitle

\begin{abstract}
\noindent  
In this paper, we give a further study on $B$-tensors and introduce doubly $B$-tensors that contain $B$-tensors.
We show that they have similar properties, including their decompositions and strong relationship with strictly (doubly) diagonally dominated tensors.
As an application, the properties of $B$-tensors are used to localize real eigenvalues of some tensors, which would be very useful in verifying the positive semi-definiteness of a tensor.

\vspace{2mm}
\noindent {\bf Key words:}\hspace{2mm} $B$-tensors, doubly $B$-tensors, decomposition of $B$-tensors, decomposition of doubly $B$-tensors, eigenvalues of tensors, positive semi-definiteness. \vspace{3mm}

\noindent {\bf AMS subject classifications (2010):}\hspace{2mm}
15A18; 15A69
\vspace{3mm}

\end{abstract}

\section{Introduction}
\hspace{12pt}
Denote $[n] := \{ 1, \ldots, n \}$.
A real $m$th order $n$-dimensional tensor (hypermatrix) $\A =
(a_{i_1\cdots i_m})$ is a multi-array of real entries $a_{i_1\cdots
i_m}$, where $i_j \in [n]$ for $j \in [m]$. Denote the set of all
real $m$th order $n$-dimensional tensors by $T_{m, n}$.  Then $T_{m,
n}$ is a linear space of dimension $n^m$.
Let $\A = (a_{i_1\cdots i_m}) \in T_{m, n}$. If the entries
$a_{i_1\cdots i_m}$ are invariant under any permutation of their
indices, then $\A$ is called a {symmetric tensor}.  Denote the
set of all real $m$th order $n$-dimensional tensors by $S_{m, n}$.
Then $S_{m, n}$ is a linear subspace of $T_{m, n}$.

For a tensor $\A =(a_{i_1 \ldots i_m}) \in T_{m, n}$ and a vector $\x = (x_1, \ldots, x_n)^\top \in \mathbb{C}^n$,
let $\A \x^{m-1}$ be a vector in $\mathbb{C}^n$ whose $i$th component is defined as
$$
(\A \x^{m-1})_i =\sum_{i_2, \ldots, i_m=1}^n a_{i i_2 \ldots i_m} x_{i_2} \ldots x_{i_m},
$$
and let $\x^{[m-1]}= (x_1^{m-1}, \ldots, x_n^{m-1})^\top$.
If there is a number $\lambda \in \mathbb{C}$ and a nonzero vector $\x \in \mathbb{C}^n$ such that
$$
\A \x^{m-1} = \lambda \x^{[m-1]},
$$
then $\lambda$ is called an eigenvalue of the tensor $\A$ and $\x$ is called an eigenvector of $\A$ corresponding to the eigenvalue $\lambda$.
We call an eigenvalue of $\A$ an H-eigenvalue of $\A$ if it has a real eigenvector $\x$.
This definition was first introduced by Qi \cite{Q}. Several other types of eigenvalues and eigenvectors for tensors were also defined in \cite{Q, L}.
In recent years, the study of tensors and the spectra of tensors (and hypergraphs) with their
various applications has attracted extensive attention and interest.

Let $\A =(a_{i_1 \ldots i_m}) \in T_{m, n}$ and $\x \in \mathbb{R}^n$. Then $\A \x^m$ is a homogeneous polynomial of degree $m$, defined by
$$
\A \x^m = \sum_{i_1, \ldots, i_m =1}^n a_{i_1 \ldots i_m} x_{i_1} \ldots x_{i_m}.
$$
Assume that $m$ is even.  If $\A \x^m \ge 0$
for all $\x \in \mathbb{R}^n$, then we say that $\A$ is {positive
semi-definite}. If $\A \x^m > 0$ for all $\x \in \mathbb{R}^n, \x \not =
\0$, then we say that $\A$ is {positive definite}. Clearly, if $m$
is odd, there is no nontrivial positive semi-definite tensors.
In \cite{Q}, Qi showed that an even order real symmetric tensor is positive semi-definite (positive definite) if and
only if all of its H-eigenvalues are nonnegative (positive).
To the best of our knowledge, positive definiteness and semi-definiteness of real symmetric tensors and their corresponding homogeneous polynomials
have applications in automatical control \cite{BM,HH,JM,Q}, magnetic resonance imaging \cite{CDHS,HHNQ,QYW,QYX} and spectral hypergraph theory \cite{HQ,LQY,Q2}.

In general, it is an NP-hard problem to verify the positive (semi-)definiteness of a real symmetric tensor \cite{HL}.
However, we can give checkable sufficient conditions for positive definite or semi-definite tensors if they have some special structure.
For instance, it is well-known that a diagonally dominated symmetric matrix is positive semi-definite.
An even order diagonally dominated symmetric tensor is positive semi-definite \cite{QS}.
It is easy to check that a given symmetric tensor is diagonally dominated or not.
This is the reason why some structured tensors are considered in the literature.

In the matrix literature, $P$-matrices are a class of important structured matrices, which were first studied systematically
by Fiedler and Pt\'{a}k \cite{FP} and found applications in linear complementarity
problems, variational inequalities and nonlinear complementarity problems and so on.
It is well-known that a symmetric $P$-matrix is positive definite \cite{CPS}.
In 2001, Pe\~{n}a proposed and studied $B$-matrices, and showed
that the class of $B$-matrices is a subclass of $P$-matrices \cite{P1}.
Thus, a symmetric $B$-matrix is positive definite.
This is another checkable sufficient condition for positive definite matrices
since it is easy to check a given matrix is a $B$-matrix or not.
Moreover, the properties of $B$-matrices are used to localize the real eigenvalues
of a real matrix and the real parts of all eigenvalues of a real matrix \cite{P1}.
On the other hand, recently, Song and Qi extended the concept of $P$-matrices and $B$-matrices to $P$-tensors,
$B$-tensors and $B_0$-tensors, obtained some nice properties about these tensors in \cite{SQ}.
By using a new approach, Qi and Song \cite{QS} showed that an even order symmetric $B$-tensor is positive definite
and an even order symmetric $B_0$-tensor is positive semi-definite.
This gives a checkable sufficient condition for positive (semi-)definite tensors.
It was also shown that these results can be extended to a general case \cite{YY}.
One natural question is: can we exploit the properties of $B$-tensors to localize the real eigenvalues of a real tensor?
This is the motivation to write this paper.
In this paper, we give a further study on $B$-tensors and by applying the properties of $B$-tensors, we derive the location of H-eigenvalues of a real symmetric tensor.
This would be very useful in verifying the positive semi-definiteness of a real symmetric tensor.
The results about $B$-tensors will be summarized in Section 2, which can be useful in the following analysis.

The concept of doubly $B$-matrices was also introduced by Pe\~{n}a \cite{P2}. It has been proved that a $B$-matrix is a doubly $B$-matrix
and a doubly $B$-matrix is a $P$-matrix. By using the properties of doubly $B$-matrices,
Pe\~{n}a also derived a similar alternative to the Brauer¡¯s theorem on ovals of Cassini \cite{P2}.
As a natural generalization, we extend the concept of doubly $B$-matrices to doubly $B$-tensors in Section 3.
It is clear that a $B$-tensor is a doubly $B$-tensor, but the converse is not true in general.
We show that a (doubly) $B$-tensor has a strong relationship with a strictly (doubly) diagonally dominated $Z$-tensor.
Similar to $B$-tensors, we also show that every principal sub-tensor of a doubly $B$-tensor is also a doubly $B$-tensor.

To characterize $B$-matrices and doubly $B$-matrices, it has been shown \cite{AT,P1} that they can be decomposed into two simple matrices in certain classes.
In Section 4, we generalize these results to the tensor case. Here,
we give the decompositions of $B$-tensors and doubly $B$-tensors, respectively.
We show that a $B$-tensor can be expressed as the sum of a $Z$-tensor and a nonnegative tensor which are both $B$-tensors, and a doubly $B$-tensor can also be expressed as the sum of a $Z$-tensor and a nonnegative tensor of simple form which are both doubly $B$-tensors.
The application of $B$-tensors is presented in Section 5.
We use the properties of $B$-tensors to localize real eigenvalues of a real tensor with some structure.
These structured tensors include $Z$-tensors, even order symmetric tensors, odd order tensors and so on.
In Section 6, we make some final remarks and raise some further questions.

Throughout this paper, we assume that $m \geq 2$ and $n \geq 1$.
We use small letters $x, y, \ldots, $ for scalars, small bold letters $\x, \y, \ldots, $ for vectors,
capital letters $A, B, \ldots, $ for matrices, calligraphic letters $\A, \B, \ldots, $ for tensors.
All the tensors discussed in this paper are real.

\section{Preliminaries}
\hspace{12pt}
Recall that a tensor $\A =(a_{i_1 \ldots i_m})\in T_{m,n}$ is called a $B$-tensor iff for all $i \in [n]$,
$$
\sum_{ i_2,\ldots,i_m=1}^n a_{i i_2 \ldots i_m} >0,
$$
and
$$
\frac{1}{n^{m-1}} \left( \sum_{ i_2,\ldots,i_m=1}^n a_{i i_2 \ldots i_m} \right) >  a_{i j_2 \ldots j_m} \text{ for all } (j_2,\ldots,j_m) \neq (i,\ldots,i).
$$
For each row $i$, denote
\begin{equation}\label{e1}
r^+_{i_\A}  = \max\{0, a_{i j_2 \ldots j_m} : (j_2,\ldots,j_m) \neq (i,\ldots,i) \},
\end{equation}
and
\begin{equation}\label{e2}
r^-_{i_\A}  = \min\{0, a_{i j_2 \ldots j_m} : (j_2,\ldots,j_m) \neq (i,\ldots,i) \}.
\end{equation}
If the content is unambiguous, we just write $r^+_i$ and $r^-_i$ for simplicity.
By definition, it is clear that $a_{i \ldots i} > r^+_i $ for all $i \in [n]$ if $\A $ is $B$-tensor.
In \cite{SQ}, it was proved that a $B$-tensor can be characterized by the following theorem.
\begin{Theorem}\label{t1}
A tensor $\A =(a_{i_1 \ldots i_m})\in T_{m,n}$ is a $B$-tensor if and only if for each $i \in [n]$,
\begin{equation}\label{e2.25}
\sum_{ i_2,\ldots,i_m=1}^n a_{i i_2 \ldots i_m} > n^{m-1} r^+_i,
\end{equation}
i.e.,
\begin{equation}\label{e2.5}
(a_{i i \ldots i} - r^+_i) > \sum_{ (i_2,\ldots,i_m )\neq (i,\ldots,i)} (r^+_i - a_{i i_2 \ldots i_m} ).
\end{equation}
\end{Theorem}

A tensor $\A =(a_{i_1 \ldots i_m})\in T_{m,n}$ is called a $Z$-tensor iff all of its off-diagonal
entries are non-positive, i.e., $a_{i_1 \ldots i_m} \leq 0$ for all $ (i_1,\ldots,i_m) \neq (i,\ldots,i)$ \cite{ZQZ};
$\A$ is called strictly diagonally dominated iff for all $i \in [n]$,
\begin{equation}\label{e2.55}
a_{i i \ldots i} > \sum \left\{ |a_{i i_2 \ldots i_m}| : (i_2,\ldots,i_m) \neq (i,\ldots,i)  \right\}.
\end{equation}
Tensor $\A$ is called strictly doubly diagonally dominated iff for all $i \in [n]$, $a_{i i \ldots i} >0$, and for all $i \neq j $ in $[n]$,
$$
a_{i i \ldots i} a_{j j \ldots j} > \sum_{(i_2,\ldots,i_m) \neq (i,\ldots,i)}  |a_{i i_2 \ldots i_m}|
\sum_{(j_2,\ldots,j_m) \neq (j,\ldots,j)}  |a_{j ij_2 \ldots j_m}|  .
$$
Note that the definition of strictly doubly dominated tensors may be different from Definition 7 in \cite{LL}.
Clearly, $\A$ is strictly doubly diagonally dominated if $\A$ is strictly diagonally dominated,
and the converse is not true in general.
It was proved in \cite{ZQZ} that a diagonally dominated $Z$-tensor is an $M$-tensor, and a strictly
diagonally dominated $Z$-tensor is a strong $M$-tensor. The definition of $M$-tensors may be
found in \cite{DQW,ZQZ}. Strong $M$-tensors are called nonsingular $M$-tensors in \cite{DQW}.

In \cite{SQ}, the relationship between these structured tensors is also given as follows.
\begin{Proposition}\label{p1}
Let $\A$ be a $Z$-tensor. Then $\A$ is a $B$-tensor if and only if $\A$ is strictly diagonally dominated.
\end{Proposition}

A tensor $\C \in T_{m,r} \  (1 \leq r \leq n)$  is called a principal sub-tensor of a tensor $\A = (a_{i_1 \cdots i_m}) \in
T_{m,n} $ iff there is a nonempty subset $J$ that composed of $r$ elements in $[n]$ such that
$$ C = (a_{i_1 \cdots i_m}), \quad \forall i_1, i_2, \cdots , i_m \in J.$$
The concept was first introduced and used in \cite{Q} for symmetric tensor. We denote by $\A^J_r$
the principal sub-tensor of a tensor $\A = (a_{i_1 \cdots i_m}) \in T_{m,n}$ such that the entries of $\A^J_r$ are indexed by
$J \subseteq [n]$ with $|J| = r \ (1 \leq r \leq n)$, and denote by $\x_J$ the $r$-dimensional sub-vector of a
vector $\x \in \mathbb{C}^n$, with the components of $\x_J$ indexed by $J$. Note that for $r = 1$, the principal
sub-tensors are just the diagonal entries.
Like the matrix case, it was also established in \cite{SQ} that all the principal sub-tensors of a $B$-tensor are also $B$-tensors.
\begin{Proposition}\label{p2}
All the principal sub-tensors of a $B$-tensor are also $B$-tensors.
\end{Proposition}

\section{Doubly $B$-tensors}
\hspace{5mm}
In this section, we give the definition of doubly $B$-tensor which is generalized from doubly $B$-matrix. The properties of this kind of structured tensor are also studied.

\begin{Definition}\label{d1}
A tensor $\A =(a_{i_1 \ldots i_m})\in T_{m,n}$ is called a doubly $B$-tensor if the following properties are satisfied:

(1) for all $i \in [n]$, $$ a_{i \ldots i} > r^+_i ,$$

(2) for all $i \neq j$ in $[n]$,
\begin{equation}\label{e2.6}
(a_{i i \ldots i} - r^+_i)(a_{j j \ldots j} - r^+_j) > \sum_{ (i_2,\ldots,i_m )\neq (i,\ldots,i)} (r^+_i - a_{i i_2 \ldots i_m} ) \sum_{ (j_2,\ldots,j_m )\neq (j,\ldots,j)} (r^+_j - a_{j j_2 \ldots j_m} ),
\end{equation}
where $r^+_i$ is defined in (\ref{e1}).
\end{Definition}

It is obvious that a $B$-tensor is a doubly $B$-tensor.
When $m=2$, a doubly $B$-tensor reduces to a doubly $B$-matrix which was proposed in \cite{P2}.
Note that another definition related to doubly $B$-tensors was proposed by Li et al. in \cite{LL}.
Besides the constraints (1) and (2) in Definition \ref{d1}, their definition \cite{LL} also assumes that for all $i \in [n]$,
$$
a_{i i \ldots i} -r^+_i \geq  \sum_{ (i_2,\ldots,i_m )\neq (i,\ldots,i)} (r^+_i - a_{i i_2 \ldots i_m} ).
$$
However, it is worth mentioning that the definition of doubly $B$-matrices proposed in \cite{P2} do not have these constraints.
So our definition can be seen as a natural generalization of the definition of doubly $B$-matrices.
Like the matrix case, it will be shown that doubly $B$-tensors share many similar properties with $B$-tensors.
On the other hand, as a general class of $B$-tensors, $MB$-tensors were introduced in \cite{LQL}.
In fact, a doubly $B$-tensor may be not an $MB$-tensor, which can be shown by the counter example $\A =(a_{i_1 i_2 i_3 i_4})\in T_{4,2}$ in \cite{LL}, namely,
$$ a_{1111}=a_{2222} =2, \quad a_{1222}=a_{2122}=a_{2212}=a_{2221} =-1 , $$
and $a_{i_1 i_2 i_3 i_4} =0$ otherwise. It is clear that $\A$ is a doubly $B$-tensor by definition, but not an $MB$-tensor since $\A$ is not positive definite. 

\begin{Proposition}\label{p3}
Suppose that $\A = (a_{i_1 \ldots i_m}) \in T_{m,n}$ is a $Z$-tensor. Then $\A$ is a doubly $B$-tensor if and only if $\A$ is strictly doubly diagonally dominated.
\end{Proposition}
\noindent{\bf Proof.}  Since $\A$ is a $Z$-tensor, we have $r^+_i = 0$ for all $i \in [n]$. By definition, the conclusion follows immediately.
\ep

\medskip
Given a tensor $\A =(a_{i_1 \ldots i_m}) \in T_{m,n}$, let $\A^+ = (b_{i_1 \ldots i_m}) \in T_{m,n} $ be the tensor defined as
\begin{equation}\label{e3}
b_{i_1 \ldots i_m} = a_{i_1 \ldots i_m} - r^+_{i_1},
\end{equation}
where $r^+_i$ is defined in (\ref{e1}). Clearly, $\A^+$ is a $Z$-tensor.

In the following, we establish a relationship between $\A$ and $\A^+$ in the case of $B$-tensors and doubly $B$-tensors.

\begin{Proposition}\label{p4}
$\A$ is a $B$-tensor if and only if $\A^+$ is a $B$-tensor.
\end{Proposition}
\noindent{\bf Proof.} Suppose that $\A =(a_{i_1 \ldots i_m}) \in T_{m,n}$ and let $\A^+ = (b_{i_1 \ldots i_m}) \in T_{m,n}$ be the tensor defined by (\ref{e3}).
Since $\A^+$ is a $Z$-tensor, we have $r^+_{i_{\A^+}} =0$ for all $i \in [n]$. Then $\A$ is a $B$-tensor if and only if for all $i \in [n]$,
$$
\sum_{ i_2,\ldots,i_m=1}^n a_{i i_2 \ldots i_m} > n^{m-1} r^+_{i_\A},
$$
i.e.,
$$
\sum_{ i_2,\ldots,i_m=1}^n b_{i i_2 \ldots i_m} = \sum_{ i_2,\ldots,i_m=1}^n \left( a_{i i_2 \ldots i_m} - r^+_{i_\A} \right) >0 = n^{m-1} r^+_{i_{\A^+}}.
$$
This is equivalent to say that $\A^+$ is a $B$-tensor.
\ep

\begin{Proposition}\label{p5}
$\A$ is a doubly $B$-tensor if and only if $\A^+$ is a doubly $B$-tensor.
\end{Proposition}
\noindent{\bf Proof.} Suppose that $\A =(a_{i_1 \ldots i_m}) \in T_{m,n}$ and let $\A^+ = (b_{i_1 \ldots i_m}) \in T_{m,n}$ be the tensor defined by (\ref{e3}).
Then $\A$ is a doubly $B$-tensor if and only if for all $i \in [n]$, $a_{i \ldots i} > r^+_{i_\A}$,  and for all $i \neq j$ in $[n]$,
$$
(a_{i i \ldots i} - r^+_{i_\A})(a_{j j \ldots j} - r^+_{j_\A}) > \sum_{ (i_2,\ldots,i_m )\neq (i,\ldots,i)} (r^+_{i_\A} - a_{i i_2 \ldots i_m} ) \sum_{ (j_2,\ldots,j_m )\neq (j,\ldots,j)} (r^+_{j_\A} - a_{j j_2 \ldots j_m} ).
$$
That is for all $i \in [n]$, $b_{i \ldots i} > 0$,  and for all $i \neq j$ in $[n]$,
$$
b_{i i \ldots i} b_{j j \ldots j}  > \sum_{ (i_2,\ldots,i_m )\neq (i,\ldots,i)} ( - b_{i i_2 \ldots i_m} ) \sum_{ (j_2,\ldots,j_m )\neq (j,\ldots,j)}
(- b_{j j_2 \ldots j_m} ).
$$
Taking into account that $\A^+$ is a $Z$-tensor, this means that $\A^+$ is a doubly $B$-tensor.
\ep

\medskip
Then, we have the following corollaries.
\begin{Corollary}\label{c1}
$\A$ is a $B$-tensor if and only if $\A^+$ is strictly diagonally dominated.
\end{Corollary}

\begin{Corollary}\label{c2}
$\A$ is a doubly $B$-tensor if and only if $\A^+$ is strictly doubly diagonally dominated.
\end{Corollary}

Similar with $B$-tensors, we show that every principal sub-tensor of a doubly $B$-tensor is also a doubly $B$-tensor.
\begin{Theorem}\label{t2}
Suppose that $\A = (a_{i_1 \ldots i_m}) \in T_{m,n}$ is a doubly $B$-tensor. Then, every principal sub-tensor of $\A$ is also a doubly $B$-tensor.
\end{Theorem}
\noindent{\bf Proof.} Let $J$ be a nonempty subset of $[n]$ with $| J | =r$ and let $\B =\A^J_r \in T_{m,r}$ be the principal sub-tensor of $\A$.
Since $\A$ is a doubly $B$-tensor, we have $a_{i \cdots i} > r^+_{i_\A} \geq r^+_{i_\B} $ for all $i \in J$.
On the other hand, for all $i \neq j$ in $J$,
\begin{eqnarray*}
(a_{i i \ldots i} - r^+_{i_\B})(a_{j j \ldots j} - r^+_{j_\B}) & \geq &  (a_{i i \ldots i} - r^+_{i_\A})(a_{j j \ldots j} - r^+_{j_\A})   \\
                                                               & >   &  \sum_{ (i_2,\ldots,i_m )\neq (i,\ldots,i)} (r^+_{i_\A} - a_{i i_2 \ldots i_m} ) \sum_{ (j_2,\ldots,j_m )\neq (j,\ldots,j)} (r^+_{j_\A} - a_{j j_2 \ldots j_m} )   \\
 & \geq & \sum_{ \begin{subarray}{c} {i_2,\ldots,i_m \in J} \\ {(i_2,\ldots,i_m )\neq (i,\ldots,i)} \end{subarray} } (r^+_{i_\A} - a_{i i_2 \ldots i_m} ) \sum_{\begin{subarray}{c} {j_2,\ldots,j_m \in J} \\ {(j_2,\ldots,j_m )\neq (j,\ldots,j)} \end{subarray} } (r^+_{j_\A} - a_{j j_2 \ldots j_m} )  \\
 & \geq & \sum_{ \begin{subarray}{c} {i_2,\ldots,i_m \in J} \\ {(i_2,\ldots,i_m )\neq (i,\ldots,i)} \end{subarray} } (r^+_{i_\B} - a_{i i_2 \ldots i_m} ) \sum_{\begin{subarray}{c} {j_2,\ldots,j_m \in J} \\ {(j_2,\ldots,j_m )\neq (j,\ldots,j)} \end{subarray} } (r^+_{j_\B} - a_{j j_2 \ldots j_m} ) .
\end{eqnarray*}
By definition, it follows that $\B$ is a doubly $B$-tensor.
\ep

\section{Decompositions of $B$-tensors and doubly $B$-tensors}
\hspace{12pt}
In \cite{QS}, Qi and Song proved that a symmetric $B$-tensor can always be decomposed to the sum
of a strictly diagonally dominated symmetric $M$-tensor and several
positive multiples of partially all one tensors. In \cite{YY}, this result was extended to a $B$-tensor
whose positive entries are invariant under any permutation.
Recently, another kind of decomposition for (doubly) $B$-matrices was introduced in \cite{AT}.
In this section, we generalize these results to the tensor case.

\medskip
Before giving the decomposition of $B$-tensors, we need the following lemma.
In fact, this result can be also derived from the fact that the set of all $B$-tensors is a convex cone.
\begin{Lemma}\label{l1}
The sum of two $B$-tensors is still a $B$-tensor.
\end{Lemma}
\noindent{\bf Proof.}
Let $\A=(a_{i_1 \ldots i_m}) \in T_{m,n}$ and $\B=(b_{i_1 \ldots i_m}) \in T_{m,n}$ be two $B$-tensors and $\C=(c_{i_1 \ldots i_m}) \in T_{m,n}$ be the sum of $\A$ and $\B$.
By Theorem \ref{t1}, we have that for all $i \in [n]$,
$$
\sum_{ i_2,\ldots,i_m=1}^n a_{i i_2 \ldots i_m} > n^{m-1} r^+_{i_\A}  \quad \text{ and } \quad \sum_{ i_2,\ldots,i_m=1}^n b_{i i_2 \ldots i_m} > n^{m-1} r^+_{i_\B}.
$$
On the other hand, it is easy to check that for all $i \in [n]$,
$$
r^+_{i_\A} + r^+_{i_\B} \geq r^+_{i_\C} \geq 0.
$$
It follows that for all $i \in [n]$,
$$
\sum_{ i_2,\ldots,i_m=1}^n c_{i i_2 \ldots i_m} = \sum_{ i_2,\ldots,i_m=1}^n a_{i i_2 \ldots i_m} + \sum_{ i_2,\ldots,i_m=1}^n b_{i i_2 \ldots i_m}
> n^{m-1} ( r^+_{i_\A} + r^+_{i_\B} ) \geq  n^{m-1} r^+_{i_\C}.
$$
This shows that $\C$ is still a $B$-tensor.
\ep

\begin{Theorem}\label{t3}
Let $\A \in T_{m,n}$. Then the following conditions are equivalent:
\begin{itemize}
\item[1.] $\A$ is a $B$-tensor.  \vspace{-2mm}
\item[2.] $\A = \B + \C$, where $\B$ is a $Z$-tensor and a $B$-tensor and $\C$ is a nonnegative $B$-tensor.
\end{itemize}
\end{Theorem}
\noindent{\bf Proof.}
1$\Rightarrow$2: Let $\I \in T_{m,n}$ be the identical tensor and let $\A^+=(b_{i_1 \ldots i_m}) \in T_{m,n}$ be the tensor defined in (\ref{e3}).
Since $\A$ is a $B$-tensor, by Corollary \ref{c1}, we can see that $\A^+$ is a strictly diagonally dominated $Z$-tensor, that is, for all $i \in [n]$,
$$
b_{ii \ldots i} >  \sum_{i_2, \ldots, i_m =1}^n (- b_{i i_2 \ldots i_m} ).
$$
It is trivial that there exists $\epsilon > 0$ such that for all $i \in [n]$,
$$
b_{ii \ldots i} - \epsilon >  \sum_{i_2, \ldots, i_m =1}^n (- b_{i i_2 \ldots i_m} ),
$$
i.e., there exists $\epsilon > 0$ such that $\A^+ - \epsilon \I $ is still a strictly diagonally dominated $Z$-tensor.
Let $\B = \A^+ - \epsilon \I$ and $\C = \A -\B $. By Proposition \ref{p1}, $\B$ is still a $B$-tensor.
It is easy to check that $\C$ is a nonnegative $B$-tensor. So the proof is completed.

2$\Rightarrow$1: The conclusion is easy to be derived according to Lemma \ref{l1}.
\ep

\medskip
In the following, we give the decomposition of doubly $B$-tensors. Unlike $B$-tensors, the sum of two doubly $B$-tensors may be not a doubly $B$-tensor.
A courter-example was given in \cite{AT} for the matrix case, see Example 2.1 of \cite{AT}.
However, we still have a similar result.
\begin{Theorem}\label{t4}
Let $\A \in T_{m,n}$. Then the following conditions are equivalent:
\begin{itemize}
\item[1.] $\A$ is a doubly $B$-tensor.  \vspace{-2mm}
\item[2.] $\A = \B + \C$, where $\B$ is a $Z$-tensor and a doubly $B$-tensor and $\C =(c_{i_1 \ldots i_m}) \in T_{m,n}$ is a nonnegative doubly $B$-tensor of the form
    \begin{equation}\label{e4.5} c_{i i_2 \ldots i_m} = \left\{ \begin{array}{ll}  c_i+\epsilon  &\   \text{ \rm if } (i_2, \dots i_m) = (i,\ldots,i), \\
     c_i  &\   \text{ \rm otherwise}, \end{array} \right.  \end{equation}
     with $c_i \geq 0 $ for $i \in [n]$ and $\epsilon > 0$.

\end{itemize}
\end{Theorem}

To prove this theorem, we need the following two simple lemmas.
\begin{Lemma}\label{l2}
Let $\A = (a_{i_1 \ldots i_m}) \in T_{m,n}$ be a doubly $B$-tensor and let $\C =(c_{i_1 \ldots i_m}) \in T_{m,n}$ be a nonnegative tensor of the form
$ c_{i i_2 \ldots i_m} = c_i $ with $c_i \geq 0 $ for $i \in [n]$.
Then, $\A +\C$ is a doubly $B$-tensor.
\end{Lemma}
\noindent{\bf Proof.}
Let $\B=\A+\C = (b_{i_1 \ldots i_m}) \in T_{m,n}$. By Proposition \ref{p5}, we need to prove that $\B^+ =(d_{i_1 \ldots i_m}) \in T_{m,n}$
is a doubly $B$-tensor. By definition, we can see that for all $i \in [n]$, $r^+_{i_\B} = r^+_{i_\A}+ c_i$.
On the other hand, for all $i, i_2, \ldots, i_m \in [n]$, we have
$$
d_{i i_2 \ldots i_m} = b_{i i_2 \ldots i_m} - r^+_{i_\B} = ( a_{i i_2 \ldots i_m} + c_i) - (r^+_{i_\A} + c_i)= a_{i i_2 \ldots i_m} - r^+_{i_\A}.
$$
This means that $\B^+ = \A^+$. Since $\A$ is a doubly $B$-tensor, by Proposition \ref{p5}, $\A^+$ is a doubly $B$-tensor. The conclusion follows immediately.
\ep

\begin{Lemma}\label{l3}
Let $\A = (a_{i_1 \ldots i_m}) \in T_{m,n}$ be a doubly $B$-tensor and let $\C =(c_{i_1 \ldots i_m}) \in T_{m,n}$ be a nonnegative diagonal tensor of the form
$ c_{i i \ldots i} = c_i $ with $c_i \geq 0 $ for $i \in [n]$.
Then, $\A + \C $ is a doubly $B$-tensor.
\end{Lemma}
\noindent{\bf Proof.}
Let $\B=\A + \C  = (b_{i_1 \ldots i_m}) \in T_{m,n}$. It is easy to see that $r^+_{i_\B} = r^+_{i_\A}$ for all $i \in [n]$.
Since $\A$ is a doubly $B$-tensor, we have for all $i \in [n]$, $a_{i i \ldots i} >r^+_{i_\A}$, and
for all $i \neq j$ in $[n]$,
$$
(a_{i i \ldots i} - r^+_{i_\A})(a_{j j \ldots j} - r^+_{j_\A}) > \sum_{ (i_2,\ldots,i_m )\neq (i,\ldots,i)} (r^+_{i_\A} - a_{i i_2 \ldots i_m} ) \sum_{ (j_2,\ldots,j_m )\neq (j,\ldots,j)} (r^+_{j_\A} - a_{j j_2 \ldots j_m} ).
$$
Hence, for all $i \in [n]$, $b_{i i \ldots i} = a_{i i \ldots i} + c_i > r^+_{i_\B}$, and
for all $i \neq j$ in $[n]$,
\begin{align*}
(b_{i i \ldots i} - r^+_{i_\B})(b_{j j \ldots j} - r^+_{j_\B}) &=     (a_{i i \ldots i} + c_i - r^+_{i_\B})(a_{j j \ldots j} + c_j- r^+_{j_\B}) \\
                                                               &\geq  (a_{i i \ldots i} - r^+_{i_\B})(a_{j j \ldots j} - r^+_{j_\B})   \\
                                                               &>
\sum_{ (i_2,\ldots,i_m )\neq (i,\ldots,i)} (r^+_{i_\B} - a_{i i_2 \ldots i_m} ) \sum_{ (j_2,\ldots,j_m )\neq (j,\ldots,j)} (r^+_{j_\B} - a_{j j_2 \ldots j_m} ) \\
&=\sum_{ (i_2,\ldots,i_m )\neq (i,\ldots,i)} (r^+_{i_\B} - b_{i i_2 \ldots i_m} ) \sum_{ (j_2,\ldots,j_m )\neq (j,\ldots,j)} (r^+_{j_\B} - b_{j j_2 \ldots j_m} ).
\end{align*}
By definition, $\B$ is a doubly $B$-tensor.
\ep

\medskip
Now we are ready to prove Theorem \ref{t4}.

\medskip
\noindent{\bf Proof of Theorem \ref{t4}.}
1$\Rightarrow$2: Let $\I \in T_{m,n}$ be the identical tensor and let $\A^+=(b_{i_1 \ldots i_m}) \in T_{m,n}$ be the tensor defined in (\ref{e3}).
Since $\A$ is a doubly $B$-tensor, by Corollary \ref{c2}, we can see that $\A^+$ is a strictly doubly diagonally dominated $Z$-tensor, that is, for all $i \in [n]$, $b_{ii \ldots i} > 0$, and for all $i \neq j$ in $[n]$,
$$
b_{ii \ldots i} b_{jj \ldots j} >  \sum_{i_2, \ldots, i_m =1}^n (- b_{i i_2 \ldots i_m} ) \sum_{j_2, \ldots, j_m =1}^n (- b_{j j_2 \ldots j_m} ).
$$
It is trivial that there exists $\delta > 0$ such that for all $i \neq j$ in $[n]$,
$$
(b_{ii \ldots i} - \delta) (b_{jj \ldots j}- \delta) > \sum_{i_2, \ldots, i_m =1}^n (- b_{i i_2 \ldots i_m} ) \sum_{j_2, \ldots, j_m =1}^n (- b_{j j_2 \ldots j_m} ).
$$
Let $\epsilon$ be chosen such that
$$
0 < \epsilon < \min\left\{ \delta, \ b_{ii \ldots i} \ | \ i\in [n] \right\}.
$$
For such $\epsilon$, we can see that $\A^+ - \epsilon \I $ is still a strictly doubly diagonally dominated $Z$-tensor.
Let $\B = \A^+ - \epsilon \I$ and $\C = \A -\B $. By Proposition \ref{p3}, $\B$ is still a doubly $B$-tensor.
It is obvious that $\C$ is a nonnegative tensor and has the form of (\ref{e4.5}) with $c_i = r^+_{i_\A}$ for $i \in [n]$.
Now we prove $\C$ is also a doubly $B$-tensor. It is easy to see that $\C^+ = \epsilon \I $, which is a doubly $B$-tensor.
By Proposition \ref{p5}, $\C$ is also a doubly $B$-tensor. So the proof is completed.

2$\Rightarrow$1: The conclusion is easy to be derived according to Lemma \ref{l2} and Lemma \ref{l3}.
\ep

\section{Application to the location of real eigenvalues }
\hspace{12pt}
In this section, the properties of $B$-tensors are applied to the location of real eigenvalues.
In order to do this, we need to define a function that acts on $T_{m,n}$.

Recall that for a tensor $\A = (a_{i_1 \ldots i_m}) \in T_{m,
n}$, the $k$th row tensor $\A_k = (a^{(k)}_{i_1 \ldots i_{m-1}}) \in T_{m-1, n}$ is
defined by $a^{(k)}_{i_1 \ldots i_{m-1}} \equiv a_{k i_1 \ldots
i_{m-1}}$, see \cite{CQ}.
The sign function sign($x$) is defined as
$$
{\rm sign}(x)= \left\{\begin{array}{cl}  1 &\   x>0, \\ 0 &\  x= 0, \\  -1 &\  x <0.  \end{array} \right.
$$

\begin{Definition}\label{d2}
Suppose that $\A = (a_{i_1 \ldots i_m}) \in T_{m,n}$ is a tensor with the $k$th row tensor $\A_k$, $k \in [n]$.
The function $F: T_{m,n} \rightarrow T_{m,n}$ is defined as
$$ F(\A)_k ={\rm sign}(a_{k \ldots k}) \A_k , \quad \forall k \in [n].$$
\end{Definition}

By definition, $F(\A)$ is a (doubly) $B$-tensor if $\A$ is a (doubly) $B$-tensor.
But the converse is not true in general. For example, if $F(\A)$ is a $B$-tensor,
the diagonal entries of $\A$ might be negative. It follows that $\A$ might be not a $B$-tensor.
In the following, we characterize the tensor $\A$ when $F(\A)$ is a (doubly) $B$-tensor, respectively.
Given a tensor $\A = (a_{i_1 \ldots i_m}) \in T_{m,n}$, let $r_i$ be defined as
\begin{equation}\label{e8}
r_i= \left\{\begin{array}{ll}  r_i^+ &\   \text{ if } a_{i \ldots i}>0, \\ 0  &\   \text{ if } a_{i \ldots i}=0, \\ r_i^- &\  \text{ if }  a_{i \ldots i}< 0,  \end{array} \right.
\end{equation}
where $r_i^+$ and $r_i^-$ are defined in (\ref{e1}) and (\ref{e2}).

\begin{Proposition}\label{p6}
Let $\A = (a_{i_1 \ldots i_m}) \in T_{m,n}$ be a real tensor and let $r_i$ be defined in (\ref{e8}).
Then $F(\A)$ is a $B$-tensor if and only if for all $i \in [n]$, $|a_{i \ldots i}| > |r_i|$ and
\begin{equation}\label{e9}
| a_{i \ldots i} -r_i|  > \sum_{ (i_2,\ldots,i_m )\neq (i,\ldots,i)} | r_i - a_{i i_2 \ldots i_m} |.
\end{equation}
\end{Proposition}
\noindent{\bf Proof.} For $i \in [n]$, by definition, the $i$th row tensor of $F(\A)$
is given by $\A_i$ if $a_{i \ldots i} >0$ and by $-\A_i$ if $a_{i \ldots i} <0$, where $\A_i$ is the $i$th
row tensor of $\A$. By Theorem \ref{t1}, $F(\A)$ is a $B$-tensor if and only if (\ref{e2.5}) holds
if $a_{i \ldots i} >0$ and
$$
 - a_{i \ldots i} - (-r_i^- ) > \sum_{ (i_2,\ldots,i_m )\neq (i,\ldots,i)} [ - r_i^- -(-a_{i i_2 \ldots i_m})]
$$
if $a_{i \ldots i} < 0$, where $r_i^+$ and $r_i^-$ are defined in (\ref{e1}) and (\ref{e2}).
Both cases are equivalent to (\ref{e9}) and the result follows.
\ep

\begin{Proposition}\label{p7}
Let $\A = (a_{i_1 \ldots i_m}) \in T_{m,n}$ be a real tensor and let $r_i$ be defined in (\ref{e8}).
Then $F(\A)$ is a doubly $B$-tensor if and only if for all $i \in [n]$, $|a_{i \ldots i}| > |r_i|$ and
for all $i \neq j$ in $[n]$,
\begin{equation}\label{e10}
| a_{i \ldots i} -r_i| | a_{j \ldots j} -r_j| > \sum_{ (i_2,\ldots,i_m )\neq (i,\ldots,i)} | r_i - a_{i i_2 \ldots i_m} | \sum_{ (j_2,\ldots,j_m )\neq (j,\ldots,j)} | r_j - a_{j j_2 \ldots j_m} |.
\end{equation}
\end{Proposition}
\noindent{\bf Proof.} Let $r_i^+$ and $r_i^-$ be defined in (\ref{e1}) and (\ref{e2}).
For $i \in [n]$, by definition, the $i$th row tensor of $F(\A)$
is given by $\A_i$ if $a_{i \ldots i} >0$ and by $-\A_i$ if $a_{i \ldots i} <0$, where $\A_i$ is the $i$th
row tensor of $\A$. We can deduce from the definition of $r_i$ that
$|a_{i \ldots i}| > |r_i|$ if and only if $a_{i \ldots i} > r_i^+$ if $a_{i \ldots i} > 0$ and $a_{i \ldots i} < r_i^-$ if $ a_{i \ldots i} < 0$.
By definition, $F(\A)$ is a doubly $B$-tensor if and only if

1) for all $i \in [n]$, $|a_{i \ldots i}| > |r_i|$;

2) if $a_{i \ldots i} , a_{j \ldots j} >0$, the inequality (\ref{e2.6}) holds;

3) if $a_{i \ldots i} >0 , a_{j \ldots j} <0$, the inequality
$$
(a_{i i \ldots i} - r^+_i)( r^-_j - a_{j j \ldots j}) > \sum_{ (i_2,\ldots,i_m )\neq (i,\ldots,i)} (r^+_i - a_{i i_2 \ldots i_m} ) \sum_{ (j_2,\ldots,j_m )\neq (j,\ldots,j)} ( a_{j j_2 \ldots j_m} - r^-_j )
$$
holds;

4) if $a_{i \ldots i} <0 , a_{j \ldots j} <0$, the inequality
$$
( r^-_i - a_{i i \ldots i} )( r^-_j - a_{j j \ldots j}) > \sum_{ (i_2,\ldots,i_m )\neq (i,\ldots,i)} (a_{i i_2 \ldots i_m} -r^-_i ) \sum_{ (j_2,\ldots,j_m )\neq (j,\ldots,j)} ( a_{j j_2 \ldots j_m} - r^-_j )
$$
holds.

We can see that all cases are equivalent to (\ref{e10}) and the result follows immediately.
\ep

\medskip
With these preparations, we are now ready to apply $B$-tensors to the location of real eigenvalues.
First, we apply $B$-tensors to the location of real eigenvalues of a $Z$-tensor.
\begin{Theorem}\label{t5}
Let $\A = (a_{i_1 \ldots i_m}) \in T_{m,n}$ be a $Z$-tensor and
let $\lambda$ be a real eigenvalue of $\A$.
Then $\lambda$ lies in the following union of closed intervals
$$
\lambda \in \bigcup_{i =1}^n  \left[  \sum_{ i_2, \ldots, i_m =1 }^n  a_{i i_2 \ldots i_m},  \
a_{i i \ldots i} -  \sum_{(i_2, \ldots, i_m) \neq (i, \ldots, i)} a_{i i_2 \ldots i_m} \right].
$$
\end{Theorem}
\noindent{\bf Proof.}
Let $\B = (b_{i_1 \ldots i_m}) \in T_{m,n}$ be a tensor defined as
$$
b_{i_1 \ldots i_m} = \begin{cases} | \lambda - a_{i_1 \ldots i_m} |  & \text{if } i_1=\ldots = i_m,  \\  a_{i_1 \ldots i_m}  &  \text{otherwise} . \end{cases}
$$
Obviously, $\B$ is also a $Z$-tensor. We claim that $\B$ is not a $B$-tensor. Otherwise, by Proposition \ref{p1}, it means $\B$ is a strictly diagonally dominated tensor, i.e., for all $i \in [n]$,
$$
| \lambda - a_{i \ldots i} | >  - \sum_{(i_2, \ldots, i_m) \neq (i, \ldots, i)} a_{i i_2 \ldots i_m}.
$$
On the other hand, let $\x$ be the eigenvector of $\A$, corresponding to $\lambda$, and let $k$ be the index such that
$|x_k| = \max \{ |x_i| : i \in [n]  \}$. Clearly, $x_k \neq 0$. Then we have
$$
(\lambda - a_{k \ldots k}) x_k^{m-1} = \sum_{(i_2, \ldots, i_m) \neq (k, \ldots, k) } a_{k i_2 \ldots i_m} x_{i_2} \ldots x_{i_m}.
$$
It follows that
$$
|\lambda - a_{k \ldots k}| \leq - \sum_{(i_2, \ldots, i_m) \neq (k, \ldots, k) }  a_{k i_2 \ldots i_m} ,
$$
which is a contradiction. Hence, the conclusion follows immediately.
\ep

\medskip
In fact, the conclusion above can also be derived by Theorem 2 of \cite{QS}.
Here, we give an alternative proof by applying $B$-tensors.
It is well known that the Laplacian tensor of a uniform hypergraph is a $Z$-tensor. For more details about Laplacian tensors and hypergraphs,
one can refer to \cite{Q2}.
As a result, the theorem above gives a lower bound and a upper bound for any real eigenvalue of the Laplacian tensor of a uniform hypergraph.

\begin{Corollary}\label{c3}
Suppose that $G$ is a $m$-uniform hypergraph with $n$ vertices. Let $\LL$ be the Laplacian tensor of $G$ and
let $d_i$ be the degree of the vertex $i$ for all $i \in [n]$.
If $\lambda$ is a real eigenvalue of $\LL$, then there exists $k \in [n]$ such that
$$
0 \leq  \lambda  \leq 2 d_k .
$$
\end{Corollary}

\medskip
In the following, we apply $B$-tensors to location of H-eigenvalues, which are real eigenvalues with real eigenvectors.
\begin{Lemma}\label{l4}
Suppose that $\A = (a_{i_1 \ldots i_m}) \in T_{m,n}$ is a $B$-tensor.
Then there does not exist a nonzero vector $\x \in \mathbb{R}^n$ such that $\A \x^{m-1} =\0$ if one of the following conditions holds:

(C1) $n=2$;

(C2) $m$ is odd;

(C3) $m$ is even and $\A$ is symmetric.
\end{Lemma}

\noindent{\bf Proof.}
Let $\A^+ = (b_{i_1 \ldots i_m}) \in T_{m,n}$ be the tensor defined by (\ref{e3}).
Suppose that $\x \in \mathbb{R}^n$ is a vector such that $\A \x^{m-1} = \0$. We need to prove $\x = \0$ under one of the conditions C1-C3.
Since $\A$ is a $B$-tensor, by Corollary \ref{c1}, $\A^+$ is a strictly diagonally dominated $Z$-tensor, which
is also a strong $M$-tensor. We may assume that $\sum_{i=1}^n x_i \neq 0$. Otherwise, we have
\begin{equation}\label{e11}
\0 = \A \x^{m-1} = \A^+ \x^{m-1} + \left( \begin{array}{c} r^+_1 \\ \vdots \\ r^+_n \end{array} \right) ( x_1 + \ldots + x_n)^{m-1}
=\A^+ \x^{m-1},
\end{equation}
where $r_i^+$ is defined in (\ref{e1}). From the proof of Theorem \ref{t5}, one can derive $\x = \0$. Hence, the conclusion follows immediately.

Since $\sum_{i=1}^n x_i \neq 0$, let $$ \y = \frac{1}{ \sum_{i=1}^n x_i } \x .$$
From (\ref{e11}), we have $\A^+ \y^{m-1} \leq \0$ and $\sum_{i=1}^n y_i = 1$.
It is equivalent to prove $\y = \0$ under one of the conditions C1-C3.
First, we prove $\y = \0$ when C1 or C2 holds.
Suppose $\y \neq \0$. Let $k$ be the index such that $$ | y_k | = \max_{i \in [n]} |y_i| .$$
Obviously, $y_k \neq 0$. Since $\A^+ \y^{m-1} \leq \0$, we have
$$
b_{k k \ldots k} y_k^{m-1} \leq \sum_{(i_2, \ldots, i_m) \neq (k, \ldots, k)} -b_{k i_2 \ldots i_m} y_{i_2} \ldots y_{i_m}.
$$
If $n=2$, we have $y_k >0$ since $\sum_{i=1}^n y_i = 1$. If $m$ is odd, we have $y_k^{m-1} >0$.
In both cases, one can derive that
$$
b_{k k \ldots k} \leq \sum_{(i_2, \ldots, i_m) \neq (k, \ldots, k)} -b_{k i_2 \ldots i_m} ,
$$
which contradicts with fact that $\A^+$ is a strictly diagonally dominated $Z$-tensor.
So $\y = \0$ when C1 or C2 holds.

Second, we prove $\x = \0$ under the condition C3. It was proved in \cite{QS} that all the H-eigenvalues of an even order symmetric $B$-tensor
are positive. This implies that $0$ is not an eigenvalue of $\A$ when C3 holds. Therefore, one can obtain that $\x = \0$ since $\A \x^{m-1} =\0$.
So the proof is completed.
\ep

\medskip
Note that unlike $B$-matrices, for a $B$-tensor $\A \in T_{m,n}$, there may exist a vector $\x \in \mathbb{R}^n$ such that $\A \x^{m-1} =\0$
if the conditions C1-C3 are not satisfied. For example,
a $B$-tensor $\A \in T_{4,3}$ is given by:
$$
\begin{array}{l}
a_{1111}=65 \text{ and } a_{1jkl} = 64 \text{ otherwise},  \\
a_{2222}=18, \ a_{2112}=15,  \text{ and } a_{2jkl} = 16 \text{ otherwise},  \\
a_{3333}=\frac{40}{3}, \ a_{3113}=11 , \text{ and } a_{3jkl} = 12 \text{ otherwise}.
\end{array}
$$
It is easy to check that for the vector $\x=(-4, 2, 3)^\top$, we have $\A \x^3 = \0 $.

\begin{Corollary}\label{c4}
Suppose that $m$ is odd or $n=2$. Let $\A = (a_{i_1 \ldots i_m}) \in T_{m,n}$. If $F(\A)$ is a $B$-tensor,
then there does not exist a nonzero vector $\x \in \mathbb{R}^n$ such that $\A \x^{m-1} =\0$.
\end{Corollary}
\noindent{\bf Proof.}
Suppose that there exists a nonzero vector $\x \in \mathbb{R}^n$ such that $\A \x^{m-1} =\0$.
It is easy to see that $F(\A) \x^{m-1} = \0$. However, by Lemma \ref{l4}, this can not happen when
$m$ is odd or $n=2$. So the conclusion follows immediately.
\ep

\medskip
Now by Lemma \ref{l4}, we give the location of H-eigenvalues of an even order symmetric tensor,
which can be useful in verifying its positive semi-definiteness.
\begin{Theorem}\label{t6}
Let $\A = (a_{i_1 \ldots i_m}) \in T_{m,n}$ be a symmetric tensor with an even order;
let $r_i^+$ and $r_i^-$ be defined in (\ref{e1}) and (\ref{e2}) and let $\lambda$ be an H-eigenvalue of $\A$.
Then $\lambda$ lies in the following union of closed intervals
$$
\lambda \in \bigcup_{i=1}^n \bigcup_{j=1}^n  \left[ a_{i i \ldots i} - r^+_i - \sum_{ (i_2,\ldots,i_m )\neq (i,\ldots,i)} (r^+_i - a_{i i_2 \ldots i_m} ),
a_{j j \ldots j} - r^-_j + \sum_{ (j_2,\ldots,j_m )\neq (j,\ldots,j)} ( a_{j j_2 \ldots j_m} - r^-_j ) \right].
$$
\end{Theorem}
\noindent{\bf Proof.} Let $\I \in T_{m,n}$ be the identity tensor. Observe that $\A - \lambda \I$ is also a symmetric tensor with an even order and
has the same off-diagonal elements as $\A$. Since $\lambda$ is an H-eigenvalue of $\A$, by Lemma \ref{l4}, we can deduce that $\A - \lambda \I$
is not a $B$-tensor. Otherwise, there does not exist a nonzero vector $\x \in \mathbb{R}^n$ such that $(\A - \lambda \I) \x^{m-1} =\0$,
which is a contradiction.
It follows that there exists an index $i \in [n]$ such that
$$
(a_{i i \ldots i} -\lambda - r^+_i) \leq \sum_{ (i_2,\ldots,i_m )\neq (i,\ldots,i)} (r^+_i - a_{i i_2 \ldots i_m} ),
$$
that is
$$
\lambda \geq a_{i i \ldots i} - r^+_i - \sum_{ (i_2,\ldots,i_m )\neq (i,\ldots,i)} (r^+_i - a_{i i_2 \ldots i_m} ).
$$

On the other hand, it is easy to see that $\lambda \I - \A$ is also a symmetric tensor with an even order and
has the opposite off-diagonal elements as $\A$. By a similar way, one can obtain that $\lambda \I - \A$ is not a $B$-tensor.
Thus, there exists an index $j \in [n]$ such that
$$
\lambda - a_{j j \ldots j} - (- r^-_j) \leq \sum_{ (j_2,\ldots,j_m )\neq (j,\ldots,j)} [ (-r^-_j) - (- a_{j j_2 \ldots j_m}) ],
$$
that is
$$
\lambda \leq a_{j j \ldots j} - r^-_j + \sum_{ (j_2,\ldots,j_m )\neq (j,\ldots,j)} ( a_{j j_2 \ldots j_m} - r^-_j ).
$$
Therefore, the conclusion holds immediately.
\ep

\medskip
Then, we have the following corollary, which was also derived from Theorem 4 of \cite{QS}.
\begin{Corollary}\label{c5}
An even order symmetric $B$-tensor is positive definite.
\end{Corollary}

\medskip
The next theorem shows that the range in Theorem \ref{t6} can be narrowed if the condition C1 or the condition C2 holds.
\begin{Theorem}\label{t7}
Let $\A = (a_{i_1 \ldots i_m}) \in T_{m,n}$ be a tensor;
let $r_i^+$ and $r_i^-$ be defined in (\ref{e1}) and (\ref{e2}) and let $\lambda$ be an H-eigenvalue of $\A$.
If $m$ is odd or $n=2$, then $\lambda$ lies in the following union of closed intervals
$$
\lambda \in \bigcup_{i=1}^n  \left[ a_{i i \ldots i} - r^+_i - \sum_{ (i_2,\ldots,i_m )\neq (i,\ldots,i)} (r^+_i - a_{i i_2 \ldots i_m} ), \
 a_{i i \ldots i} - r^-_i + \sum_{ (i_2,\ldots,i_m )\neq (i,\ldots,i)} ( a_{i i_2 \ldots i_m} - r^-_i ) \right].
$$
\end{Theorem}
\noindent{\bf Proof.} Let $\I \in T_{m,n}$ be the identity tensor. Observe that $\A - \lambda \I$
has the same off-diagonal elements as $\A$. Since $\lambda$ is an H-eigenvalue of $\A$, by Corollary \ref{c4}, we can deduce that $F(\A - \lambda \I)$
is not a $B$-tensor when $m$ is odd or $n=2$. Otherwise, there does not exist a nonzero vector $\x \in \mathbb{R}^n$ such that $(\A - \lambda \I) \x^{m-1} =\0$, which is a contradiction.
By Proposition \ref{p6}, it follows that there exists an index $i \in [n]$ such that $a_{i i \ldots i} -\lambda =0$ or
$$
a_{i i \ldots i} -\lambda>0, \qquad a_{i i \ldots i} -\lambda - r^+_i  \leq \sum_{ (i_2,\ldots,i_m )\neq (i,\ldots,i)} (r^+_i - a_{i i_2 \ldots i_m} ),
$$
or
$$
a_{i i \ldots i} -\lambda<0, \qquad r_i^- -(a_{i i \ldots i} -\lambda ) \leq \sum_{ (i_2,\ldots,i_m )\neq (i,\ldots,i)} (a_{i i_2 \ldots i_m} -r_i^-).
$$
It is equivalent that there exists an index $i$ such that
\begin{eqnarray*}
\lambda &\in & \{ a_{ii \ldots i}\} \bigcup \left[ a_{i i \ldots i} - r^+_i - \sum_{ (i_2,\ldots,i_m )\neq (i,\ldots,i)} (r^+_i - a_{i i_2 \ldots i_m} ), \
a_{ii \ldots i} \right) \bigcup \\
        &    &  \left( a_{ii \ldots i}, \ a_{i i \ldots i} - r^-_i + \sum_{ (i_2,\ldots,i_m )\neq (i,\ldots,i)} ( a_{i i_2 \ldots i_m} - r^-_i )  \right].
\end{eqnarray*}
Therefore, the conclusion holds immediately.
\ep

At last, we give an example to show that, in some cases, the intervals derived by Theorems \ref{t6} and \ref{t7} provide tight bounds
to localize the real eigenvalues, even if the bounds obtained from Gerschgorim circles are not tight.
\begin{Example}
Let $\A \in T_{m,n}$ be the all ones tensor. Clearly, It is symmetric and has real eigenvalues $0$ and $n^{m-1}$ with eigenvectors
$(1,-1,0,\ldots,0)^\top$ and $(1,1,\ldots,1)^\top$, respectively. By Theorems \ref{t6} and \ref{t7}, we derive that
all the real eigenvalues are located in the interval $[0, n^{m-1}]$, which is sharp. On the other hand,
the interval obtained from Gerschgorim circles is $[2-n^{m-1}, n^{m-1}]$.
\end{Example}

\section{Final remarks}
\hspace{12pt}
In this paper, the properties of $B$-tensors and doubly $B$-tensors are studied. It has been shown that
$B$-tensors and doubly $B$-tensors have some similar properties on their decompositions
and strong relationship with strictly (doubly) diagonally dominated tensors.
As an application, the properties of $B$-tensors are applied to the location of real eigenvalues, which can be very useful in verifying the positive semi-definiteness of a tensor.
In fact, some of these results can be extended to $B_0$-tensors. Since the proofs are similar, we omit them here.
Besides, like Lemma \ref{l4} for $B$-tensors, if we can establish a similar result for doubly $B$-tensors, the properties of doubly $B$-tensors
can also be used to localize real eigenvalues of some real tensors. These problems are worth further research.

\section*{Acknowledgement}
\hspace{12pt}
The authors are very grateful to the editor and the anonymous referees for their valuable suggestions and
constructive comments, which have considerably improved the presentation of the paper.
We would like to thank Prof. Liqun Qi for his many valuable comments and suggestions.


\begin{thebibliography}{abc99xyz}

\bibitem{AT} Ara\'{u}jo C M, Torregrosa J R. Some results on $B$-matrices and doubly $B$-matrices. Linear Algebra Appl, 2014, 459: 101-120

\bibitem{BM} Bose N K, Modaress A R. General procedure for multivariable polynomial positivity with control applications. IEEE Trans Automat Control, 1976, 21: 596-601

\bibitem{CDHS} Chen Y, Dai Y, Han D, Sun W. Positive semidefinite generalized diffusion tensor imaging via quadratic semidefinite programming.
    SIAM J Imaging Sci, 2013, 6: 1531-1552

\bibitem{CQ} Chen Z, Qi L. Circulant tensors with applications to spectral hypergraph theory and stochastic process. J Ind Manag Optim, 2016, 12(4): 1227-1247

\bibitem{CPS} Cottle R W, Pang J S, Stone R E. The Linear Complementarity Problem. Boston: Academic Press, 1992

\bibitem{DQW} Ding W, Qi L, Wei Y. $M$-tensors and nonsingular $M$-tensors. Linear Algebra Appl, 2013, 439: 3264-3278

\bibitem{FP} Fiedler M, Pt\'{a}k V. On matrices with non-positive off-diagonal elements and positive principal minors. Czechoslovak Mathematical J, 1962, 12: 382-400

\bibitem{HH} Hasan M A, Hasan A A. A procedure for the positive definiteness of forms of even-order. IEEE Trans Automat Control, 1996, 41: 615-617

\bibitem{HHNQ} Hu S, Huang Z, Ni H, Qi L. Positive definiteness of diffusion kurtosis imaging. Inverse Probl Imaging, 2012, 6: 57-75

\bibitem{HL} Hillar C J, Lim L H. Most tensor problems are NP-hard. Journal of the ACM (JACM), 2013, 60: 45

\bibitem{HQ} Hu S, Qi L. Algebraic connectivity of an even uniform hypergraph. J Comb Optim, 2012, 24: 564-579

\bibitem{JM} Jury E I, Mansour M. Positivity and nonnegativity conditions of a quartic equation and related problems. IEEE Trans Automat Control, 1981, 26: 444-451

\bibitem{LL} Li C, Li Y. Double B-tensors and quasi-double B-tensors. Linear Algebra Appl, 2015, 466: 343-356

\bibitem{LQL} Li C, Qi L, Li Y. $MB$-tensors and $MB_0$-tensors. Linear Algebra Appl, 2015, 484: 141-153

\bibitem{LQY} Li G, Qi L, Yu G. The Z-eigenvalues of a symmetric tensor and its application to spectral hypergraph theory. Numer Linear Algebra Appl, 2013, 20: 1001-1029

\bibitem{L}  Lim L H. Singular values and eigenvalues of tensors: a variational approach. Proceedings of the IEEE InternationalWorkshop on Computational Advances in Multi-Sensor Adaptive Processing, 2005, 1: 129-132

\bibitem{P1} Pe\~{n}a J M. A class of P-matrices with applications to the localization of the eigenvalues of a real matrix. SIAM J Matrix Anal Appl, 2001, 22: 1027-1037

\bibitem{P2} Pe\~{n}a J M. On an alternative to Gerschgorin circles and ovals of Cassini. Numer Math, 2003, 95: 337-345

\bibitem{Q} Qi L. Eigenvalues of a real supersymmetric tensor. J Symbolic Comput, 2005, 40: 1302-1324

\bibitem{Q2} Qi L. H$^+$-Eigenvalues of Laplacian and Signless Laplacian Tensors. Commun Math Sci, 2014, 12: 1045-1064

\bibitem{QS} Qi L, Song Y. An even order symmetric B tensor is positive definite. Linear Algebra Appl, 2014, 457: 303-312

\bibitem{QYW} Qi L, Yu G, Wu E X. Higher order positive semi-definite diffusion tensor imaging. SIAM J Imaging Sci, 2010, 3: 416-433

\bibitem{QYX} Qi L, Yu G, Xu Y. Nonnegative diffusion orientation distribution function. J Math Imaging Vision, 2013, 45: 103-113

\bibitem{SQ} Song Y, Qi L. Properties of some classes of structured tensors. J Optim Theory Appl, 2015, 165(3): 854-873

\bibitem{YY} Yuan P, You L. Some remarks on $P$, $P_0$, $B$ and $B_0$ tensors. Linear Algebra Appl, 2014, 459: 511-521

\bibitem{ZQZ} Zhang L, Qi L, Zhou G. $M$-tensors and some applications. SIAM J Matrix Anal Appl, 2014, 35: 437-452







\end{thebibliography}
\end{document}